%
%
\input amstex
\documentstyle{amsppt}
\nologo
 
\magnification=\magstep1
 
\document

\topmatter
\title
Szeg\"o kernels for certain unbounded domains in $\Bbb C^2 .$
\endtitle
\author
Friedrich Haslinger
\endauthor
 
\endtopmatter
\vskip 3 cm
\heading
1. Introduction
\endheading
\vskip 1 cm
In this paper we consider the connection between the Szeg\"o kernel of certain
unbounded domains of $\Bbb C^2$ and the Bergman
kernels of weighted spaces of entire functions of one complex
variable.
 
Let $p:\Bbb C \longrightarrow \Bbb R_+ $ denote a $\Cal
C^1$--function
and define $\Omega_p \subseteq \Bbb C^2 $ by
$$\Omega_p = \{ (z_1,z_2) \in \Bbb C^2: \Im (z_2) > p(z_1) \} .$$
Weakly pseudoconvex domains of this kind were investigated by Nagel, Rosay,
Stein and
Wainger [10,11] , where estimates for the Szeg\"o and the Bergman
kernel of the domain were made in terms of the nonisotropic
pseudometric defined in [12,13]. For the case where $p(z)=\vert z
\vert^k \ , k \in \Bbb N ,$ Greiner and Stein [5] found an explicit
expression for the Szeg\"o kernel of $\Omega_p,$ in which one can
recognize the form of the pseudometric used for the nonisotropic
estimates (see [2,8]). If $p$ is a subharmonic function,
which depends only on
the real or only on the imaginary
part of $z,$ then one can find analogous expressions and
estimates in [9].
 
Let $H^2(\partial \Omega_p)$ denote the space of all functions $f
\in L^2(\partial \Omega_p),$ which are holomorphic in $\Omega_p$ and
such that
$$\sup_{y>0}\int_{\Bbb C} \int_{\Bbb R} \vert f(z,t+ip(z)+iy) \vert^2
\, d\lambda (z) dt < \infty ,$$
where $d\lambda $ denotes the Lebesgue measure on $\Bbb C .$
We identify $\partial \Omega_p$ with $\Bbb C \times \Bbb R ,$ and
denote by $S((z,t),(w,s)), \ z,w\in \Bbb C \ , s,t\in \Bbb R ,$ the
Szeg\"o kernel of $H^2(\partial \Omega_p) .$
 
We use the tangential Cauchy--Riemann operator on
$\partial \Omega_p $
to get an expression for the Bergman kernel $K_{\tau}(z,w) $ in the
space $H_{\tau} $ of all entire functions $f$ such that
$$\int_{\Bbb C} \vert f(z) \vert^2 \exp (-2\tau p(z))\,d\lambda
(z) < \infty \ ,$$
where $\tau >0 \ ;$ in this connection we suppose that the weight
functions $p$ have a reasonable growth behavior so that the
corresponding spaces of entire functions are nontrivial, for example
if $p(z)$ is a polynomial in $\Re z$ and $\Im z .$
 
On the other hand, if one integrates the Bergman kernels with respect
to the parameter $\tau ,$ one obtains a formula for the Szeg\"o
kernel of $H^2(\partial \Omega_p) .$
 
We apply the main result for special functions $p$ to get
generalizations of results in [5,8,9].
In [7] one can find another approach to get explicit expressions
for the Szeg\"o kernel.
Finally the Bergman kernels for the spaces $H_\tau $, where $p$ is a
function of $\Re z $, are investigated, especially their asymptotic
behavior, which leads to sharp estimates and applications to
problems considered in [7] concerning a duality problem in
functional analysis.
\proclaim{Proposition 1} Let $\tau>0 .$ Then
$$K_{\tau }(z,w)=e^{\tau (p(z)+p(w))} \int_{\Bbb R} \int_{\Bbb R}
S((z,t),(w,s)) \frac{e^{i \tau (s-t)}}{p(w)-is}\,ds dt , \tag1 $$
where the integrals are to be understood in the sense of the
Plancherel theorem, i.e. in general one has only $L^2$--convergence
of the integrals.
\endproclaim
 
The fact that the above
formula (1) is not symmetric in $z$ and $w$ is due to the $L^2$--
convergence of the integrals.
 
\proclaim{Proposition 2}
$$S((z,t),(w,s))=\int_0^{\infty}K_{\tau}(z,w)
e^{-\tau (p(z)+p(w))}e^{-i\tau (s-t)}\, d\tau .\tag2 $$
\endproclaim
\vskip 1 cm
\heading
2. Proofs of Proposition 1. and 2.
\endheading
\vskip 1 cm
 
For the proof we consider the tangential Cauchy--Riemann operator
$$L=\frac{\partial}{\partial \overline{z_1}}-2i
\frac{\partial p}{\partial \overline{z_1}}(z_1)\frac{\partial}{\partial
\overline{z_2}}$$
on $\partial \Omega_p .$ Then (see [8]) $L$ is a global tangential
antiholomorphic vector field, and
$$H^2(\partial \Omega_p)=\{f\in L^2(\partial \Omega_p)
: \ L(f)=0 \  \text{as distribution}\}.$$
 
After the usual identification of $\partial \Omega_p $ with
$\Bbb C \times \Bbb R $ the tangential Cauchy--Riemann operator
has the form
$$ L=\frac{\partial}{\partial \overline z}-i\frac{\partial p}{\partial
\overline z}\frac{\partial}{\partial t} \ .$$
For a function $f \in L^2(d\lambda (z)dt)$ let $\Cal F$ denote
the Fourier transform with respect to the variable $t \in \Bbb R
\ :$
$$(\Cal F f)(z,\tau )= \int_{\Bbb R}f(z,t)e^{-it\tau }\,dt \ .$$
Then
$$ \Cal F L \Cal F^{-1}= \frac{\partial}{\partial \overline z}+
\tau \frac{\partial p}{\partial \overline z} \ .$$
$\Cal F$ and $\Cal F^{-1}$ are to be taken in the sense of the
Plancherel theorem.
 
Now let $M$ denote the multiplication operator
$$M\ :\ L^2(d\lambda (z)dt)\longrightarrow
L^2(e^{-2tp(z)}d\lambda (z)dt) $$
defined by
$$(Mf)(z,\tau )=e^{\tau p(z)}f(z,\tau ) \ ,$$
for $f\in L^2(d\lambda (z)dt) $ . Then
one has
$$\Cal F L\Cal F^{-1} = M^{-1}\frac{\partial}{\partial \overline
z}M \ .\tag3$$
 
Let $\Cal P$ denote the orthogonal projection
$$\Cal P : L^2(d\lambda (z)dt) \longrightarrow \text{Ker} L , $$
and let $P$ be the orthogonal projection
$$P : L^2(e^{-2tp(z)}d\lambda (z)dt) \longrightarrow \text{Ker}
\frac{\partial}{\partial \overline z} .$$
For fixed $\tau >0$, let $P_{\tau}$ be the orthogonal projection
$$P_{\tau} : L^2(e^{-2\tau p(z)}d\lambda (z)) \longrightarrow
\text{Ker} \frac{\partial}{\partial \overline z} .$$
Now we claim that
$$(Pf)(z,\tau )=\left\{ \aligned (P_{\tau} f_{\tau})(z) \ \ &,\tau
>0\\ 0 \ \ &,\tau \le 0\endaligned \right\} ,$$
where $f_\tau (z) = f(z,\tau ) $, for $f\in L^2(e^{-2tp(z)}
d\lambda (z)dt) .$
In order to see this it is enough to observe that a function
$f \in L^2(e^{-2tp(z)}d\lambda (z)dt)$
holomorphic with respect to the variable
$z$ has the property $f(z,t)=0$ , for almost all $t \le 0$, which is a
consequence of our assumption on the weight function $p .$
 
The next step is to show that
$$ P = M\Cal F \Cal P \Cal F^{-1} M^{-1}.\tag4$$
Denote the right side of (4) by $Q .$ We have to show that $Q^2=Q$
and that $$\text{Ker} \frac{\partial}{\partial \overline z}
\subseteq L^2(e^{-2tp(z)}d\lambda (z)dt)$$
coincides with the image of $Q .$ The first assertion follows
directly from the definition of $Q .$ For the second assertion take
a function $f \in L^2(e^{-2tp(z)}d\lambda (z)dt)$ and use (3) to
prove that
$$\frac{\partial}{\partial \overline z} Q f = M \Cal F L \Cal P
\Cal F^{-1} M^{-1} f ,$$
the last expression is zero, since $\Cal P \Cal F^{-1} M^{-1} f \in
\text{Ker} L ,$ which implies that the image of $Q$ is contained in
$\text{Ker} \frac{\partial}{\partial \overline z} .$ To prove the
opposite inclusion set $g=Q f$ for $f\in \text{Ker} \frac{\partial}
{\partial \overline z} .$ We are finish, if we can show that $Q g = f
.$ From (3) we get now
$$L \Cal F^{-1} M^{-1} f = \Cal F^{-1} M^{-1} \frac{\partial}
{\partial \overline z} f , $$
which is zero by the assumption on $f ,$ hence $\Cal F^{-1} M^{-1} f
\in \text{Ker} L$ and therefore
$$\Cal P \Cal F^{-1} M^{-1} f = \Cal F^{-1} M^{-1} f .$$
The last equality yields
$$Q g = M \Cal F \Cal P \Cal F^{-1} M^{-1} f = M \Cal F \Cal F^{-1}
M^{-1} f = f ,$$
which proves formula (4).

For a fixed $\tau >0 $ take a function $F \in L^2(e^{-2\tau
p(z)}d\lambda (z)) $ and define
$$f(z,t)= \left\{ \aligned \chi (z)F(z) \ \ &,t\ge \tau\\
                            0 \ \ &,t<\tau\endaligned  \right\} ,$$
where $\chi$ is a nonnegative, smooth function with the properties
$(\chi (z))^2=p(z),$ for $\vert z\vert \le 1$ and $\chi (z)=1,$ for
$\vert z \vert \ge 2 .$
 
Since
$$\aligned & \int_{\Bbb C} \int_{\Bbb R} \vert f(z,t)\vert^2
e^{-2tp(z)}\,dtd\lambda(z) =
 \int_{\Bbb C} \int_{\tau}^{\infty}\vert \chi (z) F(z)\vert^2
e^{-2tp(z)}\,
dtd\lambda (z) \\
& =\int_{\Bbb C}\frac{1}{2p(z)} \vert \chi (z) F(z) \vert^2
e^{-2\tau p(z)}\, d\lambda (z)
\le \text{Const.}\int_{\Bbb C} \vert F(z) \vert^2
e^{-2\tau p(z)} \, d\lambda(z) ,\endaligned $$
it
follows that $$f\in L^2(e^{-2tp(z)}d\lambda (z)dt) .$$
 
Now we use formula (4) to obtain (1): application of the operators
$M^{-1}$ and $\Cal F^{-1}$ to the function $f$ from above yields
$$\aligned \Cal F^{-1} M^{-1} f(w,t) & = \int_{\tau}^{\infty} \chi (w) F(w)
e^{t(i\sigma -p(w))}\,dt \\
& = \frac{\chi (w) F(w)e^{-\tau (p(w) -i\sigma)}}
{p(w)- i\sigma} ,\endaligned $$
 
which is a function in $L^2(d\lambda (w)d\sigma) ,$ by the properties
of the function $\chi .$
 
The next operator in (4) is now $\Cal P $ , which is the Szeg\"o
projection, hence an application of this operator can be
expressed by integration over the Szeg\"o kernel
$S((z,t),(w,\sigma )) .$ Finally we carry out the action of the operators
$\Cal F $ and $M$ and recall the properties of the operator $P$
on the left side of (4), which imply that this operator is for a
fixed $\tau $ the Bergman projection in a weighted space of
entire functions in one variable. The function $\chi $ appears on
both sides and hence cancels out.
 In this way we get formula (1).
In order to prove (2) one writes (4) in the form
$$ \Cal P = \Cal F^{-1} M^{-1} P M \Cal F ,\tag5 $$
and applies an analogous procedure as above.

\vskip 1 cm
\heading 3. Examples \endheading
\vskip 1 cm
{\bf (a)} Let $\alpha \in \Bbb R , \alpha >0 .$ We consider the function
$p(z) =\vert z \vert^{\alpha } $ and get from [6] the following expression
for the Bergman kernel $K_{\tau }(z,w) $ in the space $H_{\tau } :$
$$ K_{\tau }(z,w) = \frac{2\pi}{\alpha} \sum_{k=0}^{\infty}
(2\tau )^{2(k+1)/\alpha }\left( \  \Gamma (2(k+1)/\alpha )\right)^{-1}
z^k\overline w^k . $$
Now we apply formula (2) to this sum and get
$$ \aligned & S((z,t),(w,s)) \\
& =\frac{2\pi }{\alpha } \sum_{k=0}^{\infty}
\left( \Gamma (2(k+1)/\alpha )\right)^{-1} \
z^k\overline w^k \  2^{2(k+1)/\alpha }
\int_0^{\infty } \tau^{2(k+1)/\alpha } e^{-\tau (\vert z
\vert^{\alpha }+\vert w\vert^{\alpha })} e^{-i\tau (s-t)}
\,d\tau , \endaligned $$
evaluation of the last integral gives
$$ \Gamma \left(\frac{2(k+1)}{\alpha } +1 \right
) \ \left[ \vert z \vert^{\alpha
}+\vert w \vert^{\alpha }+ i(s-t)\right]^{-(2(k+1)/\alpha )-1}, $$
by the functional equation of the $\Gamma $--function we have
$$ \Gamma \left(\frac{2(k+1)}{\alpha } +1\right )
 =\frac{ 2(k+1)}{\alpha } \  \Gamma
(2(k+1)/\alpha ) ,$$
hence $$ S((z,t),(w,s))=\frac{2\pi }{\alpha } \sum_{k=0}^{\infty}
\frac{2(k+1)}{\alpha }\  2^{2(k+1)/\alpha } z^k\overline w^k
\  \left[ \vert z \vert^{\alpha
}+\vert w \vert^{\alpha }+ i(s-t)\right]^{-(2(k+1)/\alpha )-1} .$$
Now we set
$$ A= \frac{1}{2}( \vert z \vert^{\alpha
}+\vert w \vert^{\alpha }+ i(s-t)) $$
and carry out the summation over $k$ with the result
$$ S((z,t),(w,s))=\frac{2\pi }{\alpha^2 }A^{-1-2/\alpha }
\left( 1-\frac{z\overline w}{A^{2/\alpha }} \right)^{-2} .$$
This generalizes a result of Greiner and Stein [5], where the
same formula appears for $\alpha \in \Bbb N $ (see also [2,8]).
 
\vskip .5 cm
 
{\bf (b)} If the weight function $p$ depends
only on the real part of $z$ and satisfies
$$\int_{\Bbb R}e^{-2p(x)+2yx}\,dx <\infty ,$$
for each $y\in \Bbb R ,$ then the Bergman kernel of
$H_{\tau}$ is given by
$$K_{\tau}(z,w)=\frac{1}{2\pi}\int_{\Bbb R}\frac{\exp( \eta
(z+\overline w))}{\int_{\Bbb R}\exp (2 (r\eta -\tau p(r)))\,dr}\,
d\eta ,\tag6  $$
or
$$K_{\tau}(z,w)=\frac{\tau}{2\pi}\int_{\Bbb R}\frac{ \exp(\tau \eta
(z+\overline w))}{\int_{\Bbb R}\exp (2\tau (r\eta - p(r)))\,dr}\,
d\eta .\tag6'  $$
This follows by a modification of methods developed in [9].
To show (6) we proceed in the following way:
 
In sake of simplicity we set $\tau =1 .$
Similar to the proofs of Proposition 1 and 2 we consider the
multiplication operator
$$M_p : L^2(d\lambda (z)) \longrightarrow L^2(e^{-2p(x)}d\lambda (z))
,$$
defined by $(M_p f)(z)=e^{p(x)}f(z) \ \  , f\in L^2(d\lambda (z)) .$ Now
a computation shows that
$$\frac{\partial }{\partial \overline z} \left ( e^{p(x)}f(z) \right
) = e^{p(x)} \left ( \frac{1}{2} \frac{\partial p}{\partial x} f +
\frac{\partial f}{\partial \overline z} \right ) ,$$
which can be expressed by the operator identity
$$L(f):=\left ( M_{-p} \frac{\partial }{\partial \overline z} M_p
\right ) (f) = \frac{1}{2} \frac{\partial p}{\partial x} f +
\frac{\partial f}{\partial \overline z} .$$
Let $\Cal F $ denote the Fouriertransform with respect to $y :$
$$\Cal F f(x,\eta )=\int_{-\infty }^{\infty }f(x,y) e^{-iy\eta }\,dy
.$$
Then in the sense of distributions we have
$$ \Cal F L(f) (x,\eta )=\frac{1}{2} \left ( e^{-p(x)+\eta x}
\frac{\partial }{\partial x} \left ( e^{p(x)-\eta x}\Cal F f(x,\eta )
\right ) \right ).$$
We set $\psi (x,\eta )=e^{p(x)-\eta x} $ and define the
multiplication operator
$$\Cal M_{\psi } : L^2(d\lambda (z)) \longrightarrow
L^2(e^{-2p(x)+2yx}d\lambda (z)) $$
by $(\Cal M_{\psi }g)(x,\eta )=\psi (x,\eta )g(x,\eta ) ,$ for
$g\in L^2(d\lambda (z)) .$ Combining this with the last results we
get
$$ L=\frac{1}{2} \Cal F^{-1} \Cal M_{-\psi } \frac{\partial
}{\partial x} \Cal M_{\psi } \Cal F ,$$
and finally
$$\frac{\partial }{\partial \overline z} =
\frac{1}{2} M_p \Cal F^{-1} \Cal M_{-\psi } \frac{\partial
}{\partial x} \Cal M_{\psi } \Cal F M_{-p} .$$
In this context we consider differentiation with respect to $x$ as
an operator
$$\frac{\partial }{\partial x} : L^2(e^{-2p(x)+2yx}d\lambda (z))
\longrightarrow L^2(e^{-2p(x)+2yx}d\lambda (z)) ,$$
in the sense of distributions.
 
Further we remark that
$\text{Ker} \frac{\partial }{\partial x} $ consists of all functions
$g\in L^2(e^{-2p(x)+2yx}d\lambda (z)) ,$ which are constant in $x .$
 
By our assumption on the weight function $p$ the space
$L^2(e^{-2p(x)+2yx}dx) $ contains the constants for each $y\in \Bbb R
.$ Let $P_y $ denote the orthogonal projection of
$L^2(e^{-2p(x)+2yx}dx)$ onto the constants and $P$ the orthogonal
projection of $L^2(e^{-2p(x)+2yx}d\lambda (z))$ onto
$\text{Ker} \frac{\partial }{\partial x} .$ Then it is easily seen
that
$$( Pg )(x,y)=P_y g_y(x) ,$$
for $g \in L^2(e^{-2p(x)+2yx}d\lambda (z)) ,$ where $g_y(x)=g(x,y) .$
 
For a fixed $y\in \Bbb R$ and a function $h\in L^2(e^{-2p(x)+2yx}dx)
$ one has
$$P_y h = \frac{(h,1)}{(1,1)} 1 =
\left ( \int_{\Bbb R }e^{-2p(x)+2yx}\,dx \right )^{-1}
\int_{\Bbb R } h(x)e^{-2p(x)+2yx}\,dx .$$
Finally let $\Cal P $ denote the orthogonal projection of
$L^2(e^{-2p(x)}d\lambda (z))$ onto $H_1 =\text{Ker}
\frac{\partial }{\partial \overline z} .$
 
With the help of the above
operator identities we readily establish now
$$\Cal P = M_p \Cal F^{-1} \Cal M_{-\psi } P
\Cal M_{\psi } \Cal F M_{-p}.$$
This identity, together with the above remarks on the orthogonal
projection $P$,
implies formula (6).
\vskip 1 cm
Using (2) one gets
$$S((z,t),(w,s))=\frac{1}{2\pi}\int_0^{\infty}
\int_{\Bbb R}\frac{\tau \exp (\tau (\eta (z+\overline w)
-p(z)-p(w)-i(s-t)))}{\int_{\Bbb R}\exp (2\tau (r\eta-p(r)))\,dr
}\, d\eta \, d\tau ,$$
which is similar to an expression in [9].
 
Now we investigate the asymptotic behavior of the integral
$$\int_{\Bbb R}\exp (2\tau (r\eta -p(r)))\,dr ,\tag7 $$
which appears in formula (6), first as a function of $\eta ,$
for $|\eta | \to \infty .$
 
We restrict our attention to the case where the weight function $p$
is of the form
$$ p(r)=\frac{|r|^{\alpha }}{\alpha } , \ \alpha >1, \ r\in \Bbb R .$$
Let $p^*$ denote the Young conjugate of $p$ which is given
by
$$p^*(\eta )=\sup_{x \ge 0} \left [ x\vert \eta \vert -p(x) \right ]
= \frac{|\eta |^{\alpha^{\prime }}}{\alpha^{\prime }},\tag8 $$
where $\frac{1}{\alpha } + \frac{1}{\alpha^{\prime }}=1 .$
Note that $p^{**}=p .$
Now we can estimate the integral (7) from
above.
$$\int_{\Bbb R}\exp (2\tau (r\eta -p(r)))\,dr =
\int_0^{\infty }\exp (2\tau (r\eta - p(r)))\,dr +
\int_{-\infty }^0 \exp (2\tau (r\eta -p(r)))\,dr. $$
Let $  \lambda > 1 $. Then we have for $\eta \ge 1 $
$$\aligned  \int_0^{\infty }\exp (2\tau (r\eta - p(r)))\,dr
& \le
\int_0^{\infty }\exp (2\tau (r\eta -\lambda \eta r
+ p^*(\lambda \eta )))\,dr \\
& = \exp (2\tau (p^*(\lambda \eta ))
\int_0^{\infty } \exp(-2\tau (\lambda -1)r\eta )\,dr \\
& =\frac{\exp (2\tau p^*(\lambda \eta )}{2\tau (\lambda -1)\eta } ,
\endaligned $$
and for the second part of the integral
$$\aligned \int_{-\infty }^0 \exp (2\tau (r\eta -p(r)))\,dr & =
\int_0^{\infty } \exp (2\tau (-r\eta -p(r)))\,dr \\
& \le \int_0^{\infty } \exp (-2\tau r\eta )\,dr \\
& = \frac{1}{2\tau \eta }. \endaligned $$
For $\eta \le -1$ we estimate in the analogous way.
 
Finally for $|\eta |<1 $ we get
$$ \int_0^{\infty }\exp (2\tau (r\eta -p(r)))\,dr \le
\int_0^{\infty }\exp (2\tau (r - p(r)))\,dr ,$$
$$\aligned \int_{-\infty }^0 \exp (2\tau (r\eta -p(r)))\,dr & =
\int_0^{\infty }\exp (2\tau (-r\eta - p(r)))\,dr  \\
& \le \int_0^{\infty }\exp (2\tau (r - p(r)))\,dr . \endaligned $$
Hence for each $\eta \in \Bbb R $ we obtain
$$\int_{\Bbb R}\exp (2\tau (r\eta -p(r)))\,dr \le C(\lambda ,\tau )
\exp (2\tau p^*(\lambda \eta )) ,$$
for each $\lambda > 1$, where $C(\lambda ,\tau )>0$ is a constant depending on
$\lambda $ and $\tau .$
 
To estimate the integral in (7) from below we denote by $\mu $ the inverse
function of the derivative $ p^{\prime } $
$$\mu (\eta ) := \left (p^{\prime }\right )^{-1} (\eta ) =
|\eta |^{1/(\alpha -1)} .$$
First suppose that $\eta \ge 0 $ and observe that $ p^{\prime } $ is strictly
increasing and that the supremum in formula (8) is attained in the point
$\mu (\eta ) ,$  hence
$$\aligned \int_{\Bbb R}\exp (2\tau (r\eta -p(r)))\,dr & \ge
\int_0^{\infty }\exp (2\tau (r\eta -p(r)))\,dr \\
& \ge \exp (2\tau (\eta (\mu (\eta )+1) - p(\mu (\eta )+1))) .\endaligned $$
Next we claim that for each $\lambda $, $0<\lambda <1 ,$ the
following inequality holds
$$ 2\tau (\eta (\mu (\eta )+1) - p(\mu (\eta )+1)) \ge
2\tau (\lambda \eta \mu (\lambda \eta )-p(\mu (\lambda \eta )))-
D(\tau ,\lambda ) , \tag9 $$
for each $\eta \ge 0 ,$ where $D(\tau ,\lambda ) >0 $ is a constant
depending on $\tau $ and $\lambda .$
 
To see this we remark that
$$ \eta (\mu (\eta )+1) - p(\mu (\eta )+1) =
\eta^{\alpha /(\alpha -1)}+\eta
-1/\alpha \left ( \eta^{1/(\alpha -1)}+1 \right
)^{\alpha } ,$$
and
$$ \lambda \eta \mu (\lambda \eta )-p(\mu (\lambda \eta ))=
( 1-1/\alpha )
\lambda^{\alpha /(\alpha -1)} \eta^{\alpha /(\alpha -1)}
.$$
It suffices to show that
$$ \left ( 1-(1-1/\alpha )\lambda^{\alpha /(\alpha -1)}
\right ) \eta^{\alpha /(\alpha -1)} + \eta \ge
1/\alpha \left ( \eta^{1/(\alpha -1)}+1 \right
)^{\alpha } - D\tilde (\lambda ), $$
for each $\eta \ge 0 ,$ where $ D\tilde (\lambda )>0 $ is a
constant depending on $\lambda .$ But this follows easily from the
fact that
$$1-(1-1/\alpha )\lambda^{\alpha /(\alpha -1)}>
1/\alpha  .$$
For $\eta <0 $ we argue in a similar way.
 
On the whole we have now proved that
$$ D(\tau ,\lambda ) \exp (2\tau p^*(\eta /\lambda )) \le
\int_{\Bbb R}\exp (2\tau (r\eta -p(r)))\,dr \le
C(\lambda ,\tau )
\exp (2\tau p^*(\lambda \eta )) , \tag10 $$
for each $\eta \in \Bbb R $ and $\lambda >1 .$
 
For the conjugate function $p^* $ one obtains by the same methods
$$ D_1(\tau ,\lambda ) \exp (2\tau p(r /\lambda )) \le
\int_{\Bbb R}\exp (2\tau (r\eta -p^*(\eta )))\,d\eta  \le
C_1(\lambda ,\tau )
\exp (2\tau p(\lambda r)) , \tag11 $$
for each $ r \in \Bbb R $ and $\lambda >1 .$

The asymptotic behavior of (7) as a function of $\tau $, $\tau  \to
\infty ,$ can be derived from [1], pg. 65 :
$$\int_{\Bbb R}\exp (2\tau (r\eta -p(r)))\,dr \asymp
\left (\frac{\tau p^{\prime \prime}(\mu (\eta ))}{2\pi } \right
)^{1/2}
\exp (2\tau p^*(\eta )) .$$
 
Let
$$\exp (2\tau \wp^*(\eta )) =
\int_{\Bbb R}\exp (2\tau (r\eta -p(r)))\,dr  .$$

Then formula (6') can be written in the form
$$K_{\tau}(z,w)=\frac{\tau }{2\pi}\int_{\Bbb R} \exp \left ( 2\tau (\eta
(\frac{z+\overline w}{2})-  \wp^*(\eta))\right ) \,d\eta .
\tag12$$
 
In view of (10) and (11) this means that the Bergman kernel
$K_{\tau}(z,w) $ is in a certain sense an analytical continuation of
the original weight $\exp (2\tau p(r)) $, namely in the form
$$\exp \left (2\tau \wp (\frac{z+\overline w}{2})\right ) .$$

For $p(z)=x^2/2$ everything can be computed explicitly:
$$\int_{\Bbb R}\exp (2\tau (r\eta -r^2/2))\,dr =
(\pi /\tau )^{1/2}\exp (\tau \eta^2) ,$$
$$K_{\tau}(z,w)=\frac{\tau}{2\pi}\exp \left(\frac{\tau}{4}(z+
\overline w)^2\right)\tag13$$
and
$$S((z,t),(w,s))= \frac{1}{2\pi} \left( \frac{1}{4}(z+\overline w)^2
-\frac{1}{8}(z+\overline z)^2-\frac{1}{8}(w+\overline w)^2
-i(s-t)\right)^{-2}\tag14$$
Applying formula (1) to the expression for the Szeg\"o kernel in (14),
we arrive again at (13), now the integral with respect to $s$
converges only in $L^2$.
 
Results of this type have also been obtained by Gindikin (see [4]
or [3] ).
 
Finally we mention an estimate for the Bergman kernel, which plays an
important role in the duality problem of [7] and which,
in itself, seems to be interesting.
 
For the Bergman kernel in formula (13)
the following condition is satisfied: for each $\tau_1>\tau $
there exists $\tau_0, 0<\tau_0<\tau ,$ such that
$$\int_{\Bbb C }\int_{\Bbb C }|K_{\tau }(z,w)|^2
\exp (-2\tau_1 p(z) -2\tau_0 p(w))\,d\lambda (z)\,d\lambda (w)
<\infty .$$
This follows by a direct computation using (13).
In the general case the
integration with respect to the variable $z $ causes no problems, as
the function $z \mapsto K_{\tau }(z,w) $ belongs to the Hilbertspace $
H_{\tau_1 } ,$ for each fixed $w .$ But, afterwards,  the integration with
respect to the variable $w $ makes difficulties, because $\tau_0<
\tau .$

\newpage

{\bf Acknowledgment.} The author would like to express his sincere
thanks to A. Nagel for valuable discussions during a conference at
the M.S.R.I. in Berkeley.
 
\vskip 2 cm
 
\heading References \endheading
\vskip 1 cm
1. N.G. de Bruijn, {\it Asymptotic methods in analysis,}
North-Holland Publishing Co., Amsterdam, 1958.
 
2. K.P. Diaz, {\it The Szeg\"o kernel as a singular integral
kernel on a family of weakly pseudoconvex domains,} Trans. Amer.
Math. Soc. {\bf 304} (1987), 147--170.
 
3. B.A. Fuks, {\it Introduction to the theory of analytic
functions in several complex variables,} (in Russian)
M.,Fizmatgiz,
Moscow, 1962.
 
4. S.G. Gindikin, {\it Analytic functions in tubular regions,}
Sov. Math.--Doklady {\bf 3} (1962), 1178--1182.
 
5. P.C. Greiner and E.M. Stein, {\it On the solvability of some
differential operators of type }$\square_b ,$ Proc. Internat. Conf.,
(Cortona, Italy, 1976--1977), Scuola Norm. Sup. Pisa, Pisa,
1978, pp. 106--165.
 
6. N. Hanges, {\it Explicit formulas for the Szeg\"o kernel for
some domains in }$\Bbb C^2 ,$ J. Functional Analysis {\bf 88}
(1990), 153--165.
 
7. F. Haslinger, {\it The Bergman kernel and duality in weighted
spaces of entire functions, } preprint PAM--{\bf 310}, Berkeley,
1986.
 
8. H. Kang, $\overline{\partial }_b $--{\it equations on certain
unbounded weakly pseudoconvex domains, } Trans. Amer. Math. Soc.
{\bf 315} (1989), 389--413.
 
9. A. Nagel, {\it Vector fields and nonisotropic metrics, }
Beijing Lectures in Harmonic Analysis, E.M. Stein, Ed.,
Princeton Univ. Press, 1986.
 
10. A. Nagel, J.P. Rosay, E.M. Stein and S. Wainger, {\it
Estimates for the Bergman and Szeg\"o kernels in certain weakly
pseudoconvex domains, } Bull. Amer. Math. Soc. {\bf 18} (1988), 55--59.
 
11. A. Nagel, J.P. Rosay, E.M. Stein and S. Wainger, {\it
Estimates for the Bergman and Szeg\"o kernels in }$\Bbb C^2 $,
Ann. of Math. {\bf 129} (1989), 113--149.
 
12.  A. Nagel, E.M. Stein and S. Wainger, {\it
Boundary behavior of functions holomorphic in domains of finite type,}
Proc. Nat. Acad. Sci. U.S.A. {\bf 78} (1981), 6596--6599.

13.  A. Nagel, E.M. Stein and S. Wainger, {\it Balls
and metrics defined by vector fields I: basic properties, } Acta
Math. {\bf 155} (1985), 103--147.
 
\vskip 1.5 cm
\parindent=0cm
INSTITUT F\"UR MATHEMATIK, UNIVERSIT\"AT WIEN, STRUDLHOFGASSE 4,
\newline
A--1090 WIEN, AUSTRIA.

\enddocument